\documentclass{article}
\usepackage[letterpaper, margin=1.5in]{geometry}
\usepackage{graphicx,color}
\usepackage{multicol}
\usepackage{amsmath}
\usepackage{amsfonts}
\usepackage{amssymb}
\usepackage{booktabs}
\usepackage{algpseudocode}
\usepackage{algorithm}
\usepackage{url}
\usepackage{color}
\usepackage{tikz}
\usepackage{caption}
\usepackage{subcaption}

\usepackage{listings}
\definecolor{dkgreen}{rgb}{0,0.6,0}
\definecolor{gray}{rgb}{0.5,0.5,0.5}
\definecolor{mauve}{rgb}{0.58,0,0.82}
\newtheorem{remark}{Remark}

\newcommand{\basealpha}{\ensuremath{\alpha^{(b)}}}
\newcommand{\maxalpha}{\ensuremath{\alpha^{(max)}}}
\newcommand{\Phil}[2]{\ensuremath{\Phi^{(#1)}_{#2}}}

\newcommand{\zl}[2]{\ensuremath{z^{(#1)}_{#2}}}
\newcommand{\xl}[2]{\ensuremath{x^{(#1)}_{#2}}}
\newcommand{\yl}[2]{\ensuremath{y^{(#1)}_{#2}}}
\newcommand{\wl}[2]{\ensuremath{w^{(#1)}_{#2}}}
\newcommand{\rl}[2]{\ensuremath{r^{(#1)}_{#2}}}
\newcommand{\gl}[2]{\ensuremath{g^{(#1)}_{#2}}}
\newcommand{\el}[2]{\ensuremath{e^{(#1)}_{#2}}}

\newcommand{\x}[1]{\ensuremath{x_{(#1)}}}
\newcommand{\y}[1]{\ensuremath{y_{(#1)}}}
\newcommand{\z}[1]{\ensuremath{z_{#1}}}
\newcommand{\N}[1]{\ensuremath{N^{(#1)}}}
\newcommand{\cf}[1]{\ensuremath{m}}

\newcommand{\A}[1]{\ensuremath{A^{(#1)}}}
\newcommand{\RR}[1]{\ensuremath{R^{(#1)}}}
\newcommand{\PP}[1]{\ensuremath{P^{(#1)}}}
\newcommand{\w}[1]{\ensuremath{w^{(#1)}}}
\newcommand{\rr}[1]{\ensuremath{r^{(#1)}}}
\newcommand{\B}[1]{\ensuremath{B^{(#1)}}}
\newcommand{\vv}[1]{\ensuremath{v^{(#1)}}}
\newcommand{\e}[1]{\ensuremath{e^{(#1)}}}
\newcommand{\g}[1]{\ensuremath{g^{(#1)}}}
\newcommand{\Phill}[1]{\ensuremath{\Phi^{(#1)}}}
\newcommand{\alphal}[1]{\ensuremath{\alpha^{(#1)}}}

\newcommand{\netsize}[1]{\ensuremath{n^{[#1]}}}
\newcommand{\wsize}[1]{\ensuremath{s^{[#1]}}}

\newcommand{\globalwsize}{\ensuremath{s}}

\title{Parallelizing Over Artificial Neural Network Training Runs with Multigrid
}

\author{
    Jacob B. Schroder\thanks{Center for Applied Scientific Computing, Lawrence Livermore National Laboratory, P.O. Box 808, L-561, Livermore,CA 94551. This work was performed under the auspices of the U.S. Department of Energy by Lawrence Livermore National Laboratory under Contract DE-AC52-07NA27344.  LLNL-JRNL-736173, {\tt schroder2}\protect@{\tt llnl.gov} }
}

\begin{document}
\pagestyle{plain}
\maketitle

\begin{abstract} Artificial neural networks are a popular and effective machine
   learning technique.  Great progress has been made on speeding up the
   expensive training phase of a neural networks, leading to highly
   specialized pieces of hardware, many based on GPU-type architectures, and
   more concurrent algorithms such as synthetic gradients.  However, the
   training phase continues to be a bottleneck, where the training data must be
   processed serially over thousands of individual training runs.  This work
   considers a multigrid reduction in time (MGRIT) algorithm that is able to
   parallelize over the thousands of training runs and converge to the
   same solution as traditional serial training. MGRIT was
   originally developed to provide parallelism for time evolution problems that
   serially step through a finite number of time-steps.  This work recasts the
   training of a neural network similarly, treating neural network training as
   an evolution equation that evolves the network weights from one step to the
   next.  Thus, this work concerns distributed computing approaches for neural
   networks, but is distinct from other approaches which seek to parallelize
   only over individual training runs.  The work concludes with supporting
numerical results that demonstrate feasibility for two model problems.  \end{abstract}

\section{Introduction}
\label{sec:intro}

Artificial neural networks are a popular and effective machine learning
technique (see \cite{De2014} for a general introduction).  Great progress has
been made on speeding up the expensive training phase of an individual network,
leading to highly specialized pieces of hardware, many based on GPU-type
architectures, and more concurrent algorithms such as synthetic gradients.
However, the training phase continues to be a bottleneck, where the training
data must be processed serially over thousands of individual training runs, thus
leading to long training times.  This situation is only exacerbated by recent
computer architecture trends where clock speeds are stagnant, and future
speedups are available only through increased concurrency.  

This situation repeats itself for the time-dependent PDE simulation community, where it
has been identified as a sequential time-stepping bottleneck \cite{Fa2014}.  In
response, interest in parallel-in-time methods, that can parallelize over the
time domain, has grown considerably in the last decade.  One of the most well
known parallel-in-time methods is Parareal \cite{LiMaTu2001}, which is
equivalent to a two-level multigrid method \cite{GaVa2007}.  This work focuses
on the multilevel approach, multigrid reduction in time (MGRIT) \cite{Fa2014}.
MGRIT uses coarser representations of the problem to accelerate the solution on
the finest level.  The coarse level equations provide coarse level error
corrections to the finest level and are complemented by a fine-scale relaxation
process that corrects fine-scale errors.  The multilevel nature of MGRIT gives
the algorithm optimal parallel communication behavior and an optimal number of
operations, in contrast to a two-level scheme, where the size of the second level
limits parallelism and optimality.  For a good introduction to the history of
parallel-in-time methods, see the recent review paper \cite{Ga2015}.

\begin{figure}
   \center
   \includegraphics[width = 0.55\textwidth]{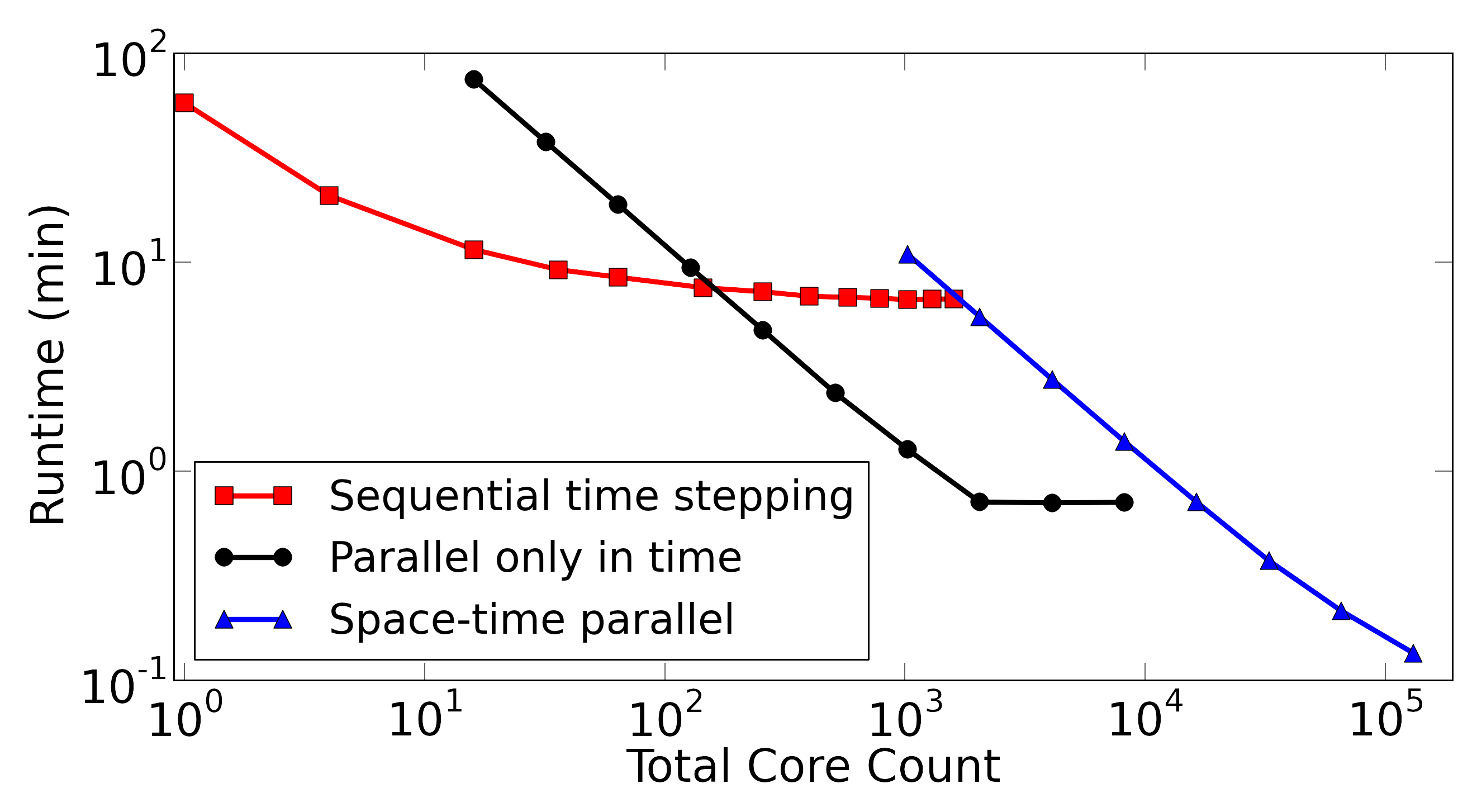}
   \caption{Wall-clock time solving a $2D$ diffusion equation for $(128)^2
   \times 16385$ space-time grid.  Sequential time-stepping is compared with
   two processor decompositions of MGRIT.}
   \label{fig:linear_example}
\end{figure}

The power of parallel-in-time methods can be seen in
Figure~\ref{fig:linear_example}.  Here, a strong scaling study for the linear
diffusion problem is carried out on the machine Vulcan (IBM BG/Q) at Lawrence
Livermore.  (See \cite{Fa2016}.)   Sequential time-stepping is compared with
two processor decompositions of MGRIT, a time-only parallel MGRIT run and a
space-time parallel MGRIT run using 64 processors in space and all other
processors in time. The maximum possible speedup is approximately 50, and the
cross-over point is about 128 processors, after which MGRIT provides a speedup. 

Two important properties of MGRIT are the following.  One, MGRIT is able to
parallelize over thousands of training runs and \emph{also} converge to the
exact same solution as traditional sequential training.  MGRIT does this by
recasting the training of a neural network as an evolution process that evolves
the network weights from one step to the next, in a process analogous to
time-steps.  Thus, this work concerns distributed computing approaches for
neural networks, but is distinct from other approaches which seek to
parallelize only over individual training runs.  The second property is that
MGRIT is nonintrusive.  The evolution operator is user-defined, and is often an
existing sequential code that is wrapped within the MGRIT framework.  MGRIT is
not completely agnostic to the user's problem; however, in that the convergence
of the method is problem dependent.  Nonetheless, the approach developed here
could be easily applied to other existing implementations of serial neural
network training.

Parallelizing over the time dimension and even over the steps of a general
evolutionary process can seem surprising at first, but the feasibility of this
approach has been shown.
For instance, MGRIT has already been successfully used in the non-PDE setting for powergrid
problems \cite{LeFaWoTo2016}.  In fact, it can be natural to think of neural
network training as a time-dependent process.  For instance, a network is
trained by ``viewing" consecutive sets of images.  In the biological world,
this would certainly be happening in a time-dependent way, and the
computer-based approximations to these networks can also be thought of as
stepping through time, as training runs are processed.  More precisely, if the
training of a network with $\globalwsize$ weights is represented by a function
application $\Phi$, traditional serial training of the network can be described as 
\begin{equation} 
   \label{eqn:naivephi}
      w_{i+1} =  \Phi(w_i), \quad \text{for some initial } w_0, \text{ and } i = 0,1,\dots,N,
\end{equation}
where $w_{i}$ is a vector of the $\globalwsize$ weight values at training step
$i$, $w_{i+1}$ represents the \emph{evolution} of these weights to step $i+1$,
and $N$ is the total number of training runs.  Other researchers have made
similar an analogies between a training run and a time-step \cite{Ho1982, Li1987}. 

Thus, the goal here is to parallelize over these $N$ consecutive $\Phi$
applications, as MGRIT has already parallelized over time-steps in other settings, e.g.
for Navier-Stokes \cite{FalKatz2014} and nonlinear parabolic problems
\cite{Fa2016}. The resulting MGRIT algorithm will converge to the \emph{same
solution} as that obtained from a traditional sequential processing of the
training runs. 

\begin{figure}[t!]
\centering
   \begin{subfigure}[t]{1.0\textwidth}
        \centering
        \includegraphics[width = 0.31\textwidth]{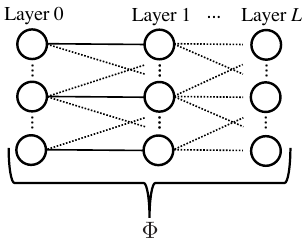}
        \caption{Visual representation of an arbitrary multilayer $\Phi$ application.}
         \label{fig:single_phi}
    \end{subfigure}\\
   \begin{subfigure}[t]{0.48\textwidth}
        \centering
        \includegraphics[width = 0.666666\textwidth]{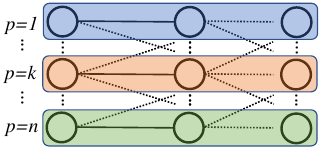}
        \caption{Parallelism inside of $\Phi$, decomposing across $n$ nodes.}
         \label{fig:lat_parallelism}
    \end{subfigure}%
    $\quad$
    \begin{subfigure}[t]{0.48\textwidth}
        \centering
        \includegraphics[width = 0.5833333\textwidth]{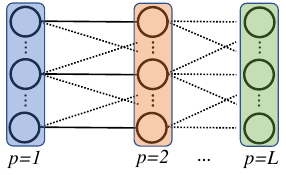}
        \caption{Parallelism inside of $\Phi$, decomposing across $L$ layers.}
         \label{fig:vert_parallelism}
    \end{subfigure}\\
    \begin{subfigure}[t]{1.0\textwidth}
         \centering
         \includegraphics[width = 0.98\textwidth]{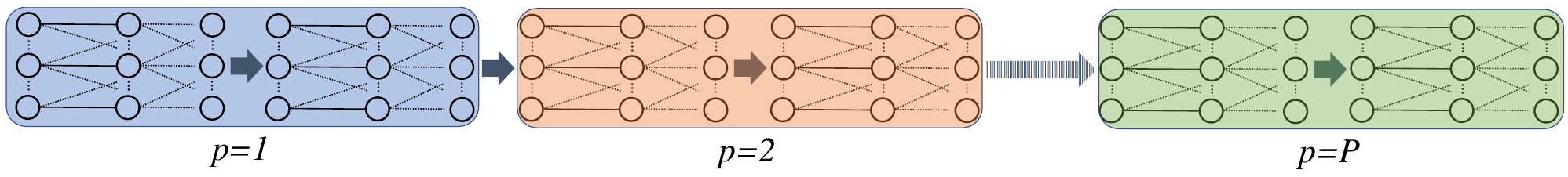}
         \caption{Parallelism across $\Phi$ applications, i.e., across training runs, for $P$ processors.} 
         \label{fig:temp_parallelism}
     \end{subfigure}
    \caption{Approaches to parallelism for artificial neural networks.}
    \label{fig:newparallelism}
\end{figure}

Figure \ref{fig:newparallelism} depicts the proposed form of parallelism over
$\Phi$ applications, and compares it to other types of possible parallelism.
Figure \ref{fig:single_phi} provides a reference figure for a single $\Phi$
application.  Figure \ref{fig:lat_parallelism} depicts parallelism inside of a
single $\Phi$ application, where the problem is decomposed across processors
(denoted ``$p$"), and each processor owns some slice of the $n$ nodes.  (While
the number of nodes is usually not constant across layers, we ignore that here
simply for the purpose of having a straight-forward pictorial representation).
Figure \ref{fig:vert_parallelism} depicts parallelism inside of a single $\Phi$
application, where the $L$ layers are decomposed across the processors.
Finally, Figure \ref{fig:temp_parallelism} depicts parallelism where entire
$\Phi$ applications are decomposed across processors, which is the approach
pursued here.  Each processor owns two $\Phi$ applications, but this
is not fixed.\footnote{These types of parallelism could be mixed, as has been
done for PDE-based problems in \cite{Fa2014,FaMaScSo2016}, but the scope of
this paper is on approach in Figure \ref{fig:temp_parallelism}.}

As motivation for MGRIT, consider a simple thought experiment for
a two-level MGRIT approach to training.  Let the fine level contain 100
training runs, and the coarse level contain 50 training runs.  MGRIT will use
the coarse level only to \emph{accelerate} convergence on the fine level by
computing error corrections on the coarse level.  It is reasonable to consider
the possibility that a 50 training run example could provide a good error
correction to a 100 training run example.  MGRIT will do this in a multilevel
way, so that the 50 training run example is never solved exactly, thus avoiding
sequential bottlenecks.

The challenge of developing strategies for good MGRIT convergence for our model
neural network problems allows this paper to make some novel achievements.  To
the best of our knowledge, this is the first successful application of
parallel-in-time techniques to neural networks.  To do this, we develop some
novel MGRIT strategies such using the learning rate in a way similar to
time-step sizes.  We also explore serializing over all
training examples, to allow for more potential parallelism.

In Section \ref{sec:motivation}, we introduce our first model neural network,
the MGRIT algorithm, and some motivation for designing an MGRIT solver for
neural network training.  In Section \ref{sec:method}, we develop the proposed
MGRIT approach, apply it to two model problems, and compute maximum possible
speedup estimates.  In Section \ref{sec:conclusion}, concluding remarks are
made.

\section{Motivation}
\label{sec:motivation}

In this section, we first formalize our notation for MGRIT and the example
neural networks, and then conclude with motivation for the proposed method. 

\subsection{Model Neural Network}
\label{sec:ann_intro}

To define a multilayer network, we let there be $\netsize{j}$ nodes on level
$j$ and $\wsize{j}$ weights connecting level $j$ to level $j+1$. In general, a
superscript $^{[\cdot]}$ will denote a level in the neural network, not MGRIT
level. We say that there are in total $\globalwsize$ weights in the network.
Indexing always starts at 0.

For the model problems considered here, the networks are fully connected that
is each node on level $j$ has a weight connecting it to every node on level
$j+1$.  Thus, we skip any further notation defining the edges (connections) in
the network.  While we assume a fixed network topology for each example, we
note that adaptive network topologies (e.g., adding and removing weights) could
be integrated into our solver framework (see \cite{FaMaScSo2016}) by using
strategies researched for MGRIT in the classical PDE setting.  

Our simplest neural network model problem is based on the three-layer example from \cite{Trask_3Layer}.
This network is trained to learned an exclusive OR (XOR) operation, where the $K$th training input is
$\x{k}$ and the $K$th expected output is $\y{k}$, i.e.,
\begin{align}
   \label{eqn:3layer_data}
   \x{0} &= \begin{bmatrix} 0 & 0 & 1 \end{bmatrix}, \quad \y{0} = \begin{bmatrix} 0 \end{bmatrix} \\
   \x{1} &= \begin{bmatrix} 0 & 1 & 1 \end{bmatrix}, \quad \y{1} = \begin{bmatrix} 1 \end{bmatrix} \nonumber\\
   \x{2} &= \begin{bmatrix} 1 & 0 & 1 \end{bmatrix}, \quad \y{2} = \begin{bmatrix} 1 \end{bmatrix} \nonumber\\
   \x{3} &= \begin{bmatrix} 1 & 1 & 1 \end{bmatrix}, \quad \y{3} = \begin{bmatrix} 0 \end{bmatrix} \nonumber.
\end{align}
The single hidden layer has four nodes, yielding an input layer of
$\netsize{0} = 3$, the hidden layer of $\netsize{1} = 4$, and a single output
node of $\netsize{2} = 1$.  The total number of weights is $\globalwsize = 16$ 
($3 \cdot 4 + 4 \cdot 1$).  The threshold function is a standard sigmoid,
\begin{equation}
   \text{threshold}(v) = \frac{v}{1 + e^{-v}},
\end{equation}
and standard backpropagation \cite{RuHiWi1985, Le1985, Pa1985} is used to train the network.  The strategy is that
MGRIT must first be effective for such simple model problems, before it is reasonable
to experiment with larger problems. 

The serial training code for this network is from \cite{Trask_3Layer} and
reproduced in Appendix \ref{sec:3layer_code} for convenience.  However, we
briefly review the general process of network training for readability. First,
the forward propagation evaluates the network by initializing the input nodes
to values from the $i$th training instance $\x{i}$.  Then these node values are
used in a matrix-vector multiply with the $\netsize{1} \times \netsize{0}$
matrix representing the weights connecting layers 0 and 1.  Then the threshold
function is applied element-wise to the result vector of size $\netsize{1}$.
This procedure of matrix-vector multiply and threshold evaluation continues
across all the hidden layers, until values at the output nodes are obtained.
The initial gradient is then obtained by differencing $\y{i}$ with the value(s)
at the output layer.  This gradient is then propagated backwards, updating the
weights between each level.  The size of the step taken in the gradient
direction is called the \emph{learning rate} $\alpha$.

The proposed MGRIT algorithm will provide the same solution as this sequential
training process, to withing some tolerance.  Later, we will apply MGRIT to a
larger, four-layer problem. 

\subsection{MGRIT}
\label{sec:mgrit_intro}

The description of MGRIT also requires some notation and nomenclature.  Let the
$\Phi$ from equation (\ref{eqn:naivephi}) be the propagator or evolution
operator that trains the network.  One \emph{training run} is a single $\Phi$
application.  Each training run is functionally equivalent to a time-step in
the traditional MGRIT framework.  Thus, we will call a training run a
``training step" from now on, to emphasize the connection between
time-dependent evolution problems and the training of artificial neural
networks.

The $i$th training step has a corresponding training data set called $\z{i}$.
For the above example in Section \ref{sec:ann_intro}, we define 
\begin{equation}
   \label{eqn:xyz}
   \z{i} = \begin{bmatrix} \x{0}, \x{1}, \x{2}, \x{3}, \y{0}, \y{1},\y{2},\y{3} \end{bmatrix} 
            \in \mathbb{R}^{16}.
\end{equation}
Thus, $\z{i}$ is the same for all $i$ here, and the network does a batch
training on all four training instances simultaneously (although we will
explore variations on this later).

In conclusion, one $\Phi$ application is equivalent to a forward propagation of
the training data in $\z{i}$ with $w_i$, followed by a backpropagation step
that creates the updated set of weights $w_{i+1}$.  When done sequentially over all
training steps, equation (\ref{eqn:naivephi}) is obtained.

However, the serial training process can also be represented as a global
operator with $N$ block rows, representing the $N$ training steps.  To do this
we ``flatten" each set of weights so that $w_i$ becomes a vector of length
$\globalwsize$.  This global operator of size $(N\globalwsize) \times (N\globalwsize)$ is
\begin{equation}
   A (w) \equiv 
      \left[ \begin{array}{cccc} I \\-\Phi() & I \\ & \ddots &  \ddots \\ & & -\Phi() & I \end{array} \right] 
      \left[ \begin{array}{c} w_0 \\ w_1 \\ \vdots \\ w_{\N{0}} \end{array} \right] =  
         \left[ \begin{array}{c} w_0 \\ 0 \\ \vdots \\ 0 \end{array} \right] \equiv g
      \label{block_tri},
\end{equation}
where $w_0$ equals the initial values assigned to the weights by the user.  A
forward solve using forward substitution of this equation is equivalent to the
traditional serial training of a neural network.  Instead, MGRIT solves
equation (\ref{block_tri}) iteratively, with multigrid reduction, allowing for
parallelism over the training steps.  That is, all training steps are computed
simultaneously.  Being an iterative method, the proposed approach also
converges to the exact same solution as that produced by traditional serial
training.  MGRIT accomplishes this by constructing a hierarchy of training
levels.  Each coarse level provides an error correction to the next finer level in
order to accelerate convergence on that level.  Complementing these error corrections
is a local relaxation process that resolves fine-scale behavior.

Regarding complexity, both serial propagation and MGRIT are $O(N)$ methods, but
MGRIT is highly parallel.\footnote{MGRIT is $O(N)$ when it converges in a
scalable way, with a fixed iteration count regardless of problem size.} On the
other hand, the computational constant for MGRIT is higher, generally leading
to a cross-over point, where beyond a certain problem size and processor count,
MGRIT offers a speedup (see Figure \ref{fig:linear_example}).

\subsection{MGRIT Algorithm}

We will first describe MGRIT for linear problems, and then describe the
straight-forward extension to nonlinear problems.  The description is for a
two-level algorithm, but the multilevel variant is easily obtained through
recursion.  

We denote $N = \N{0}$ as the number of training steps on level 0, which is the
finest level that defines the problem to be solved, and $\N{\ell}$ as the
problem size on level $\ell$.  In general, a superscript $^{(\cdot)}$ will
hereafter denote a level in the MGRIT hierarchy, and a subscript will denote a
training step number.  If the subscript denoting a training step number is
omitted, the quantity is assumed to be constant over all training steps on that
level.

\subsubsection{MGRIT for Linear Problems}

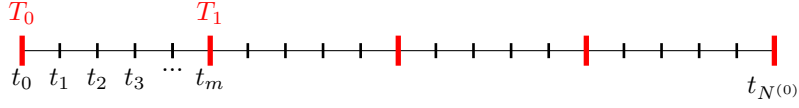
\begin{figure}
   \center
   \usetikzlibrary{decorations.pathreplacing}
   \begin{tikzpicture}
      \path[draw] (0,0) -- (10,0);

      \foreach \x in {0,1,...,3}
         \path[draw, line width = 1pt] (\x/2,0.1) -- (\x/2,-0.1) node [below] {$t_{\x}$};
         \path[draw, line width = 1pt] (2,0.1) -- (2,-0.1) node [below] {$...$};
         \path[draw, line width = 1pt] (2.5,0.1) -- (2.5,-0.1) node [below] {$t_m$};
         \foreach \x in {4,5,...,20}
            \path[draw, line width = 1pt] (\x/2,0.1) -- (\x/2,-0.1);

            \path[draw, line width = 2pt, color = red] (0,-0.2) -- (0,0.2) node [above] {$T_0$};
            \path[draw, line width = 2pt, color = red] (2.5,-0.2) -- (2.5,0.2) node [above] {$T_1$};
            \path[draw, line width = 2pt, color = red] (5,-0.2) -- (5,0.2) ;
            \path[draw, line width = 2pt, color = red] (7.5,-0.2) -- (7.5,0.2) ;
            \path[draw, line width = 2pt, color = red] (10,0.2) -- (10,-0.2) node [below] {$\textcolor{black}{t_{\N{0}}}$};
   \end{tikzpicture}
      \caption{Fine- and coarse-levels. Fine-level points (black) are present on only the fine-level, whereas coarse-level points (red) are on both the fine- and coarse-level. 
      \label{grid}}
\end{figure}

MGRIT solves an evolution problem over a finite number of evolution steps.  For
time-dependent PDEs, these steps correspond to time-steps.  For a more general
evolution problem, these steps need not be PDE-based.  Regardless, MGRIT defines
a fine level and coarse level of evolution steps, based on the coarsening factor
$\cf{0}$. (Here, we use a level independent coarsening
factor.) There are $\N{0}$ steps on level 0 and $\N{1} = \N{0}/\cf{0}$ for
level 1.  Two such levels for $\cf{0} = 4$ are depicted in 
Figure~\ref{grid}.\footnote{For time-dependent PDEs, these levels represent a fine and coarse time
step spacing, respectively, with the fine time-step size four times smaller
than the coarse level time-step size.}  Typical values for $\cf{\ell}$ are 2, 4,
8, and 16. 

Multigrid reduction is similar to an approximate block cyclic reduction
strategy, and similarly eliminates block rows and columns in the system
(\ref{block_tri}).  If all the block rows corresponding to the black fine
points in Figure \ref{grid} are eliminated from the system (\ref{block_tri}),
then using the level notation (level 1 is coarse, and level 0 is fine), the
coarse system becomes
\begin{equation}
   \label{eqn:idealcoarse}
   \A{1} \w{1} = \left [ \begin{array}{cccc} I \\-(\Phill{0})^{\cf{0}} & I \\ & \ddots &  \\ & & -(\Phill{0})^{\cf{0}} & I \end{array} \right ]\left [ \begin{array}{c} \wl{1}{0} \\ \wl{1}{1}\\ \vdots \\ \wl{1}{\N{1}} \end{array} \right ]  = \RR{0} \g{0} = \g{1},
\end{equation}
where $(\Phill{0})^{\cf{0}}$ represents the $\cf{0}$th power, i.e., $\cf{0}$
consecutive applications of \Phill{0}.  Here, we assume that $\Phill{\ell}$ is
the same function over all training steps for simplicity, although the method
extends easily to the case where $\Phill{\ell}$ depends on the training step
number.  The operator $\RR{0}$ represents ``Ideal" restriction and we also
define an interpolation $\PP{0}$, as in
\cite{Fa2014},
\begin{subequations}
   \label{eqn:RPideal}
   \begin{align}
      \RR{0} &= \begin{bmatrix} I \\  & (\Phill{0})^{\cf{0}-1} & \dots & \Phill{0} & I \\ &&&&& \ddots  \\ &&&&&& (\Phill{0})^{\cf{0}-1}  & \dots & (\Phill{0}) & I
      \end{bmatrix}, \\
      (\PP{0})^T &= \begin{bmatrix} I & (\Phill{0})^T & \dots & (\Phill{0})^{\cf{0}-1,T} \\ &&& \ddots  \\ &&&& I & (\Phill{0})^T  & \dots & (\Phill{0})^{\cf{0}-1,T} 
      \end{bmatrix} .
   \end{align}
\end{subequations}
The operator $\PP{0}$ injects coarse points from the coarse- to fine-level, and then
extends those values to fine-points with a $\Phill{0}$ application.  This will be
equivalent to injection, followed by F-relaxation (defined below).

The operator $\A{1}$ can be formed in a typical multigrid fashion with $\A{1} =
\RR{0} \A{0} \PP{0}$, and we refer to this as the ideal coarse-level operator
because the solution to (\ref{eqn:idealcoarse}) yields the exact fine-level
solution at the coarse points.\footnote{The ideal coarse-level operator can
also be obtained with the injection operator $\RR{0}_I$, $\RR{0} A \PP{0} =
\RR{0}_I A \PP{0}$. For efficiency, we therefore will use injection later on in
the algorithm.} The interpolation $P$ is also ideal, and uses $\Phi$ to compute
the exact solution at the fine points.  However, this exact solution process is
in general as expensive as solving the original fine-level problem, because of
the $(\Phill{0})^{\cf{0}}$ applications in $\A{1}$.  MGRIT ``fixes" this issue by
approximating $\A{1}$ with $\B{1}$, i.e.,
\begin{equation}
   \B{1} = \left [ \begin{array}{cccc} I \\-\Phill{1} & I \\ & \ddots &  \\ & & -\Phill{1} & I \end{array} \right ],
\end{equation}
where $\Phill{1}$ is an approximation of $(\Phill{0})^{\cf{0}}$.  One obvious
choice here is to approximate $(\Phill{0})^{\cf{0}}$ with $\Phill{0}$, i.e.,
approximate $m$ training steps with one training step.  For time-dependent PDEs,
$\Phill{1}$ is usually a single time-step using a larger time-step size that
approximates $m$ finer time-steps.  How well $\Phill{1}$ approximates
$(\Phill{0})^{\cf{0}}$ will govern how fast MGRIT converges, and choosing a
$\Phill{1}$ that allows for fast convergence is typically the main topic of
research for a given problem area.  

\begin{figure}
   \center
   \usetikzlibrary{shapes}
   \usetikzlibrary{plotmarks}

   \begin{tikzpicture}
      \draw (0,0)--(8,0);
      \filldraw (9,0) circle (0pt) node[right] {F-relaxation};
      \filldraw (-0.15,-0.15) rectangle (0.15,0.15);
      \filldraw (3.85,-0.15) rectangle (4.15,0.15);
      \filldraw (7.85,-0.15) rectangle (8.15,0.15);

      \foreach \x in {1,2,3,5,6,7}
          \filldraw (\x,0) circle (3pt) node [below,font=\large] {$\underset{g}{\uparrow}$};
          \foreach \x in {0,1,2,5,6,4}
            \draw[->, thick] (\x,0.2) to[out=90,in=90] node [above, midway] {$\Phi$} (\x+0.92,0.2); 
   \end{tikzpicture}

   \begin{tikzpicture}
      \draw (0,0)--(8,0);
      \filldraw (9,0) circle (0pt) node[right] {C-relaxation};
      \filldraw (-0.15,-0.15) rectangle (0.15,0.15);
      \filldraw (3.85,-0.15) rectangle (4.15,0.15);
      \filldraw (7.85,-0.15) rectangle (8.15,0.15);
      \foreach \x in {4,8}
          \filldraw (\x,0) circle (3pt) node [below,font=\large] {$\underset{g}{\uparrow}$};
          \foreach \x in {1,2,3,7,5,6}
              \filldraw (\x,0) circle (3pt);
              \foreach \x in {3,7}
                 \draw[->, thick] (\x,0.2) to[out=90,in=90] node [above, midway] {$\Phi$} (\x+0.92,0.2); 
   \end{tikzpicture}
   \caption{F- and C-relaxation for coarsening by factor of 4 \label{FCRelax}}
\end{figure}
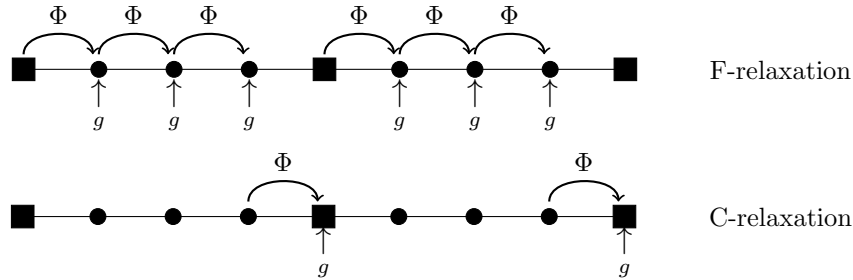

Thus, the two-level algorithm uses the coarse level to compute an error
correction based on the fine-level residual, $\rr{0} = \g{0} - \A{0}(\w{0})$,
mapped to the coarse level with $\RR{0}$
(see Step \ref{line:resid} in Algorithm~\ref{alg:mgrit}).  
Complementing the coarse level correction is a
local relaxation process on the fine-level which resolves fine-scale behavior.
Figure~\ref{FCRelax} (from \cite{Fa2016}) depicts F- and C-relaxation on a level
with $\cf{0}=4$. F-relaxation updates the solution at each F-point by
propagating the solution forward from the nearest coarse point (C-point).  And,
C-relaxation updates every C-point with a solution propagation from the nearest
F-point.  Overall, relaxation is highly parallel, with each interval of
F-points being independently updated during F-relaxation, and likewise for
C-relaxation. The nonintrusiveness of MGRIT is also apparent because the
internals of $\Phill{\ell}$ are opaque to MGRIT.  Thus, this
approach can be easily applied by wrapping existing code bases. 

\subsubsection{Nonlinear Algorithm} 
\label{sec:mgrit_nonlinear} 

The above two-level process can be easily to nonlinear problems using full
approximation storage (FAS), a nonlinear multigrid scheme \cite{Brandt_1977,
FalKatz2014}.\footnote{We mention that the F-relaxation two-level MGRIT
algorithm is equivalent to the popular parallel-in-time Parareal algorithm
\cite{GaVa2007}.}  The key addition of FAS to the multigrid process is how
the coarse-level equation is Step 3 is computed.

The two-level FAS version of MGRIT is described in Algorithm~\ref{alg:mgrit}
(taken from \cite{FalKatz2014}).  A recursive application of the algorithm at
Step 4 makes the method multilevel.  This can be done recursively
until a trivially sized level is reached, which is then
solved sequentially.  Ideal interpolation (\ref{eqn:RPideal}) is represented as
two steps (5 and 6). Initial values for the weights $\wl{0}{0}$ can in principle be
anything, even random values, which is the choice made here for $\wl{0}{0}$ as the
initial condition.  For all other points, $\wl{0}{i} = 0$ is the initial value.

\begin{algorithm}
   \caption{MGRIT($\A{0}$, $\w{0}$, $\g{0}$, $tol$)} 
\label{alg:mgrit}
\begin{algorithmic}[1]
   \State Apply F- or FCF-relaxation to $\A{0} (\w{0}) = \g{0}$.
   \State Inject the fine level approximation and its residual to the coarse level:\newline
   $
   \wl{1}{i} \leftarrow \wl{0}{mi},\quad
   \rl{1}{i} \leftarrow \gl{0}{mi} - \left( \A{0} (\w{0}) \right)_{mi}.
   $ \label{line:resid}
   \State Solve $\B{1} ( \vv{1}) = \B{1} ( \w{1} ) + \rr{1}$.
   \State Compute the coarse level error approximation: $\e{1} \simeq \vv{1} - \w{1}$.
   \State Correct $\w{0}$ at $C$-points: $\wl{0}{mi} = \wl{0}{mi} + \el{1}{i}.$
   \If{$\|\g{0} - \left( \A{0} (\w{0}) \right)\| \le tol$} \label{line:convg}
      \State Update $F$-points with apply F-relaxation to $\A{0} (\w{0}) = \g{0}$, and halt.
   \Else{}
      \State Go to step 1.
   \EndIf
\end{algorithmic}
\end{algorithm}

\begin{figure}
  \center
  \includegraphics[width = 0.6\textwidth]{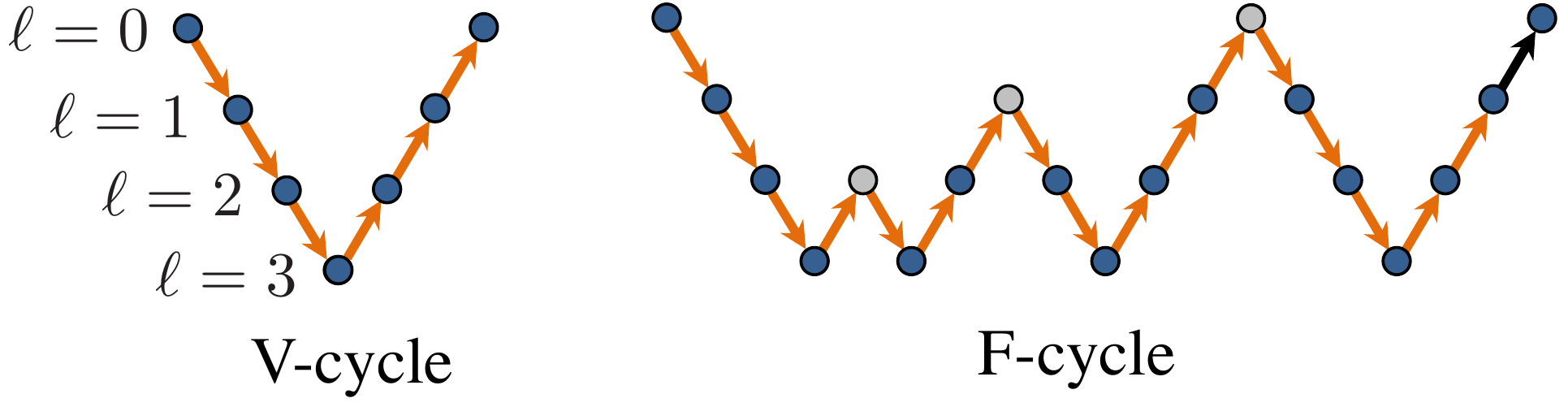}
  \caption{Example multigrid V-cycle and F-cycle with 4 levels. } 
  \label{fig:Vcycle}
\end{figure}

Multigrid solves can use a variety of cycling strategies.  Here, we consider V-
and F-cycles, as depicted in Figure~\ref{fig:Vcycle}, which shows the order in
which levels are traversed.  The algorithm continues to cycle until the convergence
check based on the fine level residual is satisfied in line \ref{line:convg}.

\subsection{Motivation}
\label{sec:sub_motivation}

The motivation for the proposed approach is based on similarities between neural network
training and time-dependent evolution problems.  We have already shown that the
sequential neural network training approach of equation (\ref{eqn:naivephi})
can represented as a global system (\ref{block_tri}) solvable by MGRIT.
However, the similarity goes even further, because we can represent any $\Phi$
operator so that it resembles a forward Euler time-step.  
\begin{equation}
   \label{eqn:fwdEuler}
   w_{i+1} = \Phi(w_i) = w_i + \alpha \hat{\Phi}(w_i),
\end{equation}
where $\hat{\Phi}$ represents the
computation of the gradient, i.e., the search direction for this update, and
$\alpha$ represents the length of the search step.  (See lines 45--47 in
Appendix \ref{sec:3layer_code}.)  But, the $\alpha$ parameter also resembles a
time-step size, which could be varied in a level-dependent manner inside of
MGRIT. This will be a key insight required for good MGRIT performance.

The other motivation is based on the fact that both the forward and backward
phases of neural network training resemble discretized hyperbolic 
PDEs, which typically have an upper or lower triangular matrix form.  If
the neural network training algorithm were a purely linear process, then the forward
phase would look like a block lower-triangular matrix-vector multiply, where
the $k$th block-row is the size of the $k$th layer.  And the backward phase, would
be analogously  upper block-triangular.  Thus, we will consider
solver strategies from \cite{Do2016} which considered MGRIT for hyperbolic
PDEs.  These strategies are slow coarsening (i.e., small $\cf{\ell}$ values),
FCF-relaxation (as apposed to F-relaxation), and F-cycles.

\section{Method}
\label{sec:method}

In this section, we formalize how we describe an MGRIT solver for artificial neural networks.
We must define $\Phil{\ell}{i}$, or the propagator on
level $\ell$ for training step $i$.  This corresponds to block row $i$ in
equation (\ref{block_tri}).  We first define a ``Naive Solver", because it is the
simplest, out-of-the-box option.  Next, we will explore modifications that
lead to faster convergence.

\subsection{Naive Solver}
\label{sec:naive}
To define $\Phil{\ell}{i}$, three chief subcomponents must be defined.
For the \textbf{Naive Solver}, we choose the following.

\paragraph{Coarsening factor $\mathbf{\cf{\ell}}$:}

The coarsening factor $\cf{\ell}$, defines the problem size on each level.  For
example, $\cf{0} = 2$, implies that $\N{1} = \N{0} / 2$.  For the Naive Solver,
we coarsen by a factor of 2. 
This yields a sequence of problems sizes on each level, $\N{\ell}$ of 
\begin{align*}
   [\N{\ell}] := [\N{0}, \N{0}/2, \N{0}/4, \N{0}/8, \N{0}/16 \dots].
\end{align*}

\paragraph{Learning rate $\mathbf{\alphal{\ell}}$:}

The learning rate on level $\ell$ is defined by $\alphal{\ell}$.  For the Naive
Solver, we choose a fixed learning rate over all levels defined by a
base learning rate $\basealpha$, i.e.,
\begin{align}
    \label{eqn:naive_alpha}
    [\alphal{\ell}] := [\basealpha, \basealpha, \basealpha, \dots].
\end{align}

\paragraph{Training set $\mathbf{\zl{\ell}{i}}$:}
The training set at level $\ell$ and training step $i$ is $\zl{\ell}{i}$.  For the examples considered here,
\begin{equation}
    \label{eqn:zldefn}
   \zl{\ell}{i} = \begin{bmatrix} \xl{\ell}{i},\; \yl{\ell}{i} \end{bmatrix},
\end{equation}
where $\xl{\ell}{i}$ is the input for the network, and $\yl{\ell}{i}$ is the desired output. 
For the Naive Solver $\zl{\ell}{i}$ will contain \emph{all} of the training
data, and is fixed, not changing relative to $\ell$ or $i$.  The data set
$\xl{\ell}{i}$ and \yl{\ell}{i} will contain four rows, and this training data
will be processed in a standard ``batch" manner for neural networks (see
Appendix \ref{sec:3layer_code}).  In other words, each application of
$\Phil{\ell}{i}$ will train over all of the training data.  However, we will
consider other options later, such as serializing over the training data.

Thus, we can say that $\Phil{\ell}{i}$ is defined by 
\begin{equation}
   \label{eqn:phi_defn}
    \Phil{\ell}{i} := \begin{cases}
                         \text{Choice for } \zl{\ell}{i}\\
                         \text{Choice for } \alphal{\ell} \\
                         \text{Choice for } \cf{\ell} \\ 
                      \end{cases}
\end{equation}
and that 
\begin{equation}
   \label{eqn:phil}
   \Phil{\ell}{i}(\wl{\ell}{i}) \rightarrow \wl{\ell}{i+1}
\end{equation}
where $\wl{\ell}{i}, \wl{\ell}{i+1} \in \mathbb{R}^{\globalwsize}$ represent
vectors of the $\globalwsize$ weights defining the network.

\subsection{Results}

For our results we will experiment using FCF-relaxation and compare 
V-cycles, F-cycles and two-level solvers, while experimenting with various
$\Phil{\ell}{i}$ definitions.  The maximum coarsest level size is 10
training steps, and the halting tolerance is $10^{-9} \sqrt{\N{0}}$.  We scale
the halting tolerance in a manner similar to PDE applications, so that the
tolerance grows commensurately with problem size.  The maximum number of
iterations is 50, thus results that took longer are denoted as ``50+" with a
corresponding average convergence rate of ``*".  The average convergence rate
$\rho$ is computed as a geometric average, $\rho = (||r_k||/||r_0||)^{1/k}$.
The quantity $r_k = g - A(w)$ is the residual from equation (\ref{block_tri})
after the $k$ MGRIT iterations, and the norm $||r_k||$ is the standard
Euclidean norm.

\paragraph{Experimental Goal:}  The goal of these experiments is to examine
various MGRIT solvers, using heuristics to choose $\Phil{\ell}{i}$.  A
successful MGRIT solver will converge quickly, in say 10 or 15 iterations, and
show little to no growth in iteration counts as the problem size increases.
(Remember, MGRIT is an iterative algorithm, which updates the solution until
the convergence check based on the fine level residual is satisfied in line
\ref{line:convg} of Algorithm \ref{alg:mgrit}.) Such a result will be said to
be ``scalable" or ``nearly scalable", because the operation count required to
solve equation (\ref{block_tri}) will be growing linearly with problem size.
In parallel, however, there will always be a log-growth in communication
related to the number of levels in the solver.  For a good introduction to
multigrid and the notion of scalability, see \cite{BrHeMc2000}.

While two-level solvers are impractical, we will first examine them (as is
typical for multigrid methods research).  If success cannot be obtained in this
simplified setting, there is little hope that a multilevel version of the
solver will perform better.

\subsection{Naive Solver Results}

Using the Naive Solver from Section \ref{sec:naive}, in a two-level setting,
yields the results in Table \ref{tab:2levelnaive}.  The learning rate is
$\basealpha = 1.0$, taken from from \cite{Trask_3Layer}.  While convergence is
poor, this result shows that MGRIT can converge for the problem in question.
\label{sec:naiveresults}
\begin{table}[!ht]
  \centering
  \begin{tabular}{c | ccccccccc }
  \toprule 
     $\N{0}$       & 100  &   200 &   400 &   800 &   1600 &   3200 &   6400 &  12800 &  25600   \\
  \midrule                                                                                    
    Iters       &  14   &  19   &  24   &  31   &  40    &  49    &  50+   &  50+   &  50+\\
    $\rho$      & 0.31  & 0.43  & 0.48  & 0.59  & 0.66   & 0.72   &   *    &   *    &   * \\
  \bottomrule
  \end{tabular}
  \caption{Three-layer model problem, Naive solver, Two-level results, $\basealpha = 1.0$}
   \label{tab:2levelnaive}
\end{table}

\subsection{Improvements to Naive}

\subsubsection{Increasing $\alphal{\ell}$}
\label{sec:incalpha}

In this section, we consider the motivation in Section
\ref{sec:sub_motivation}, and consider increasing $\alphal{\ell}$ on coarser
levels.  In particular, we consider doubling $\alphal{1}$ so that $\alphal{1} =
2 \alphal{0}$.  Remarkably, this single change to the Naive Solver yields fast,
scalable two-level iteration counts.
\begin{table}[!ht]
  \centering
  \begin{tabular}{c | ccccccccc }
  \toprule 
    $\N{0}$       & 100  & 200  & 400  & 800  & 1600 & 3200 & 6400 & 12800 & 25600 \\
  \midrule             
    Iters       &  6   &  6   &  6   &  6   &  6   &  6   &  6   &   6   &   5   \\
    $\rho$      & 0.10 & 0.10 & 0.10 & 0.10 & 0.09 & 0.09 & 0.08 &  0.08 &  0.07 \\
  \bottomrule
  \end{tabular}
  \caption{Three-layer model problem,  Solver 1, Two-level results, $\basealpha=1.0$}
  \label{tab:2levelIncreasingAlpha}
\end{table}

More generally, we consider defining $\alphal{\ell}$ as 
\begin{equation}
   \label{eqn:inc_alpha}
   \alphal{\ell} = \min(2^{\ell} \basealpha, \maxalpha),
\end{equation}
where $\maxalpha$ is some experimentally determined maximum stable
$\alphal{\ell}$ value.\footnote{This connection between the size of $\alphal{\ell}$ and stability
(i.e., $\w{\ell}$ values that do not diverge over multiple $\Phill{\ell}$
applications) shows another similarity between $\alphal{\ell}$ and time-step
size.}
Here, $\maxalpha = 8.0$, i.e., $\alphal{\ell}$ is allowed to increase up to
level $3$.  Beyond that point, the training becomes ``unstable", and
$\alphal{\ell}$ can no longer increase.  While the halt to increasing
$\alphal{\ell}$ values may negatively impact the algorithm, this halt occurs
only on coarser levels, and so its negative impact on MGRIT convergence will
hopefully be somewhat mitigated.

We denote this as \textbf{Solver 1}, which uses equation (\ref{eqn:inc_alpha})
to define $\alphal{\ell}$, $m=2$ and equation (\ref{eqn:zldefn})
to define $\zl{\ell}{i}$ as containing \emph{all} the training data.  
Together with $\basealpha=1.0$ and $\maxalpha=8.0$, this fully defines the method
explored in this section.

Table \ref{tab:FcycleIncreasingAlpha} depicts the results for Solver 1.
Initially the results look scalable (i.e.,
flat iteration counts and stable convergence rates $\rho$). But, eventually a
log-type growth in the iterations becomes apparent.  A possible culprit is the
fact that $\alphal{\ell}$ is capped and cannot continually increase.  V-cycle
results are quite poor and hence omitted for brevity.
\begin{table}[!ht]
  \centering
  \begin{tabular}{c | ccccccccc }
  \toprule 
    $\N{0}$       & 100  & 200  & 400  & 800  & 1600 & 3200 & 6400 & 12800 & 25600 \\
  \midrule             
    Iters       &  6   &  6   &  7   &  7   &  7   &  8   &  11  &   16  &   21  \\
    $\rho$      & 0.10 & 0.10 & 0.09 & 0.10 & 0.09 & 0.07 & 0.17 &  0.24 &  0.34 \\
  \bottomrule
  \end{tabular}
   \caption{Three-layer model problem, Solver 1, F-cycle results, $\basealpha=1.0$}
   \label{tab:FcycleIncreasingAlpha}
\end{table}

One important parameter that affects the MGRIT convergence rate is the base
learning rate $\basealpha$.  Decreasing this value improves the MGRIT
convergence rate, while also requiring more training steps to adequately train
the network.  Exploration of this trade-off is future research, but we make note
of it here.  If we take Solver 1, decrease $\basealpha$ from 1.0 to
0.5 and keep $\maxalpha=8.0$, we get the improved F-cycle iterations in Table
\ref{tab:FcycleIncreasingAlpha_SmallBaseAlpha}
\begin{table}[!ht]
  \centering
  \begin{tabular}{c | ccccccccc }
  \toprule 
    $\N{0}$       & 100  & 200  & 400  & 800  & 1600 & 3200 & 6400 & 12800 & 25600 \\
  \midrule                                                                       
    Iters       &  5   &  5   &  5   &  6   &  6   &  5   &  6   &   9   &   12  \\
    $\rho$      & 0.07 & 0.08 & 0.08 & 0.06 & 0.05 & 0.04 & 0.05 &  0.07 &  0.18 \\
  \bottomrule
  \end{tabular}
   \caption{Three-layer model problem, Solver 1, F-cycle results, $\basealpha=0.5$}
  \label{tab:FcycleIncreasingAlpha_SmallBaseAlpha}
\end{table}

\subsubsection{Serializing the Training Data}
\label{sec:3layer-serialize}

In this section, we note that there is unutilized parallelism in the
method. All four training examples are processed simultaneously as a ``batch".
That is, $\zl{\ell}{i}$ is defined by an $\xl{\ell}{i}$ and $\yl{\ell}{i}$ with
four rows, representing the four training sets (see equation
(\ref{eqn:zldefn}).  Instead, each of the four training sets can be broken up
and processed individually, resulting in four times as many possible steps to
parallelize.  More precisely, for $K$ training instances, define $\xl{\ell}{i}$
and $\yl{\ell}{i}$ to correspond to the mod$(i,K)$-th training instance.

Here $K=4$, so $\Phil{0}{0}$ would train the network with the
first training data set with $\zl{0}{0} = [\x{0}, \y{0}]$.  $\Phil{0}{1}$ would
train the network for the second training data set with $\zl{0}{1} = [\x{1},
\y{1}]$, and so on.  Once the training step index passes the value $K=4$, the
process repeats, so that $\Phil{0}{4}$ again trains on the first training data
set with $\zl{0}{4} = [\x{0},\y{0}]$.  This strategy exposes significantly more
parallelism than the ``unserialized" $\Phil{\ell}{i}$.  

We denote this as \textbf{Solver 2}, which defines
\begin{equation}
   \label{eqn:zl_serialized}
   \zl{\ell}{i} = [\x{j},\; \y{j}], \;\; j = \text{mod}(i,K),
\end{equation}
Solver 2 continues using equation (\ref{eqn:inc_alpha}) to define
$\alphal{\ell}$ and $m=2$.  Together with $\basealpha=1.0$ and $\maxalpha=8.0$,
this fully defines the method explored in this section.

Table \ref{tab:2levelSerializedIncreasingAlpha} depicts two-level results using Solver 2,
which show a scalable convergence rate, but the rate is slower than that in Table
\ref{tab:2levelIncreasingAlpha} for Solver 1.
\begin{table}[!ht]
  \centering
  \begin{tabular}{c | ccccccccc }
  \toprule 
    $\N{0}$       & 100  & 200  & 400  & 800  & 1600 & 3200 & 6400 & 12800 & 25600 \\
  \midrule       
    Iters       &  10  &  12  &  15  &  18  &  20  &  20  &  20  &   20  &   19  \\
    $\rho$      & 0.18 & 0.23 & 0.28 & 0.36 & 0.40 & 0.39 & 0.39 &  0.38 &  0.38 \\
  \bottomrule
  \end{tabular}
   \caption{Three-layer model problem, Solver 2, two-level results, $\basealpha=1.0$}
   \label{tab:2levelSerializedIncreasingAlpha}
\end{table}

However, Table \ref{tab:FcycleSerializedIncreasingAlpha} depicts Solver 2 results for
F-cycles, and this shows an unscalable increase in iteration counts.
\begin{table}[!ht]
  \centering
  \begin{tabular}{c | ccccccccc }
  \toprule 
    $\N{0}$     & 100  & 200  & 400  & 800  & 1600 & 3200 & 6400 & 12800 & 25600 \\
  \midrule       
    Iters       &  10  &  12  &  15  &  21  &  26  &  27  &  29  &   41  &   50  \\
    $\rho$      & 0.18 & 0.23 & 0.29 & 0.40 & 0.47 & 0.47 & 0.50 &  0.60 &  0.99 \\
  \bottomrule
  \end{tabular}
   \caption{Three-layer model problem, Solver 2, F-cycle results, $\basealpha=1.0$}
   \label{tab:FcycleSerializedIncreasingAlpha}
\end{table}

One possible contributing factor the unscalable iteration growth in Table
\ref{tab:FcycleSerializedIncreasingAlpha} is the fact that 
$\alphal{\ell}$ is allowed to only grow on the first four levels.  To investigate this
possibility, Table \ref{tab:FcycleSerializedIncreasingAlpha_SmallBaseAlpha}
runs the same experiment as in Table \ref{tab:FcycleSerializedIncreasingAlpha},
but with with $\basealpha=0.1$ and $\maxalpha=30.0$.  We see that this
combination of a smaller base $\basealpha$, which allows $\alphal{\ell}$ to grow up to level 
$\ell=8$ dramatically improves the scalability of the F-cycles.  The downside
is the smaller learning rate.  We leave a further investigation of this for
future work, because a deeper theoretical understanding is required to
meaningfully guide our tests.
\begin{table}[!ht]
  \centering
  \begin{tabular}{c | ccccccccc }
  \toprule 
    $\N{0}$     & 100  & 200  & 400  & 800  & 1600 & 3200 & 6400 & 12800 & 25600 \\
  \midrule      
    Iters       &  5   &  5   &  5   &  6   &  7   &  8   &  10  &   13  &   14  \\
    $\rho$      & 0.05 & 0.06 & 0.07 & 0.11 & 0.10 & 0.14 & 0.15 &  0.25 &  0.24 \\
  \bottomrule
  \end{tabular}
   \caption{Three-layer model problem, Solver 2, F-cycle results, $\basealpha=0.1$}
   \label{tab:FcycleSerializedIncreasingAlpha_SmallBaseAlpha}
\end{table}

\begin{remark} The reader will note that this strategy of serializing the
   training data is independent of $\ell$.  For example, MGRIT interpolation
   and restriction connects point $i$ on level 1 to point $2i$ on level 0.
   However, this serialized training set definition defines the training
   instance on level 1 for point $i$ as mod$(i,K)$ and on level 0 for point
   $2i$ as mod$(2i,K)$.  The values mod$(i,K) \neq$ mod$(2i,K)$, in general.
   Experiments were run where the training instances on level 0 were injected
   to all coarser levels, to avoid this discrepancy, but this approach lead to
   worse convergence.  It is hypothesized that this is the case because the
   complete elimination of training instances on the coarser levels leads to
   poor spectral equivalence between coarse and fine $\Phi$.  To address this,
   experiments were also run where coarse training sets represented averages of
   fine-level training sets, and this improved convergence, but not enough to
   compete with the above approach, even in a two-level setting.  The averaging
   was done such that $\xl{1}{i}$ was the average of $\xl{0}{2i}$ and
   $\xl{0}{2i+1}$, i.e., the average of the corresponding C- and F-point on
level 0.  \end{remark}

\subsubsection{Potential Speedup}
\label{sec:3layer_speedup}

While the proposed method does allow for parallelism over training steps, we
need to quantify possible parallel speedups for the proposed approach.  To do
that, we use a simple computational model that counts the \emph{dominant
sequential} component of the MGRIT process and compares
that to the dominant sequential component of the traditional serial training
approach.  The dominant sequential component in both cases is the required
number of consecutive $\Phil{\ell}{i}$ applications, because this is the main
computational kernel. 

For sequential time stepping, the number of consecutive $\Phil{0}{i}$
applications is simply $\N{0}$.  However for MGRIT, we must count the number of
consecutive $\Phil{\ell}{i}$ applications as the method traverses levels.
Consider a single level during a V-cycle with a coarsening factor of 2.
FCF-relaxation contains 3 sequential $\Phil{\ell}{i}$ applications, the sweeps
over F- then C- then F-points.  
Additionally, each level must compute an FAS residual (see Algorithm
\ref{alg:mgrit}), which requires two matrix-vector products.  In parallel at
the strong scaling limit, each process in time will own two block rows that
will likely be processed sequentially, resulting in 4 total sequential
$\Phil{\ell}{i}$ applications for the computation of the two matrix-vector products.  Last,
on all levels, but level 0, there is an additional $\Phil{\ell}{i}$ application
to fill in values at F-points for the interpolation.  The coarsest level is an
exception and only carries out a sequential solve no longer than the maximum
allowed coarsest level size (here 10). 

Taken together, this results in 7 sequential $\Phil{\ell}{i}$ applications on
level 0, and 8 on all coarser levels, except the coarsest level which has
approximately 10.  Thus, the dominant sequential component for a V-cycle is
modeled as 
\begin{equation}
   \label{eqn:seq_partV}
   \sigma_V(L) = \text{niter} \left(17 + \sum_{\ell = 1}^{L-2} 8\right),
\end{equation}
where $\text{niter}$ is the number of MGRIT iterations, or V-cycles, required for
convergence.

The dominant sequential component of an F-cycle is essentially the concatenation 
of successive V-cycles starting on ever finer levels (see Figure \ref{fig:Vcycle}), i.e.,
\begin{equation}
   \label{eqn:seq_partF}
   \sigma_F(L) = \text{niter} \sum_{\ell = 1}^{L-2} \sigma_V(\ell).
\end{equation}
The reader will note that there is no final V-cycle done on level 0 in equation
(\ref{eqn:seq_partF}).

Thus, we define the potential speedup for V-cycles as $\N{0} / \sigma_V$ and
$\N{0} / \sigma_F$ for F-cycles.  This is a useful speedup measure, but it
measures speedup relative to an equivalent sequential training run.  For
example, for the case of the double $\Phill{\ell}$ applications, this measure
will compute a speedup against a sequential run that also does a double
$\Phill{\ell}$ application.

\begin{figure} \center \includegraphics[width =
   0.6\textwidth]{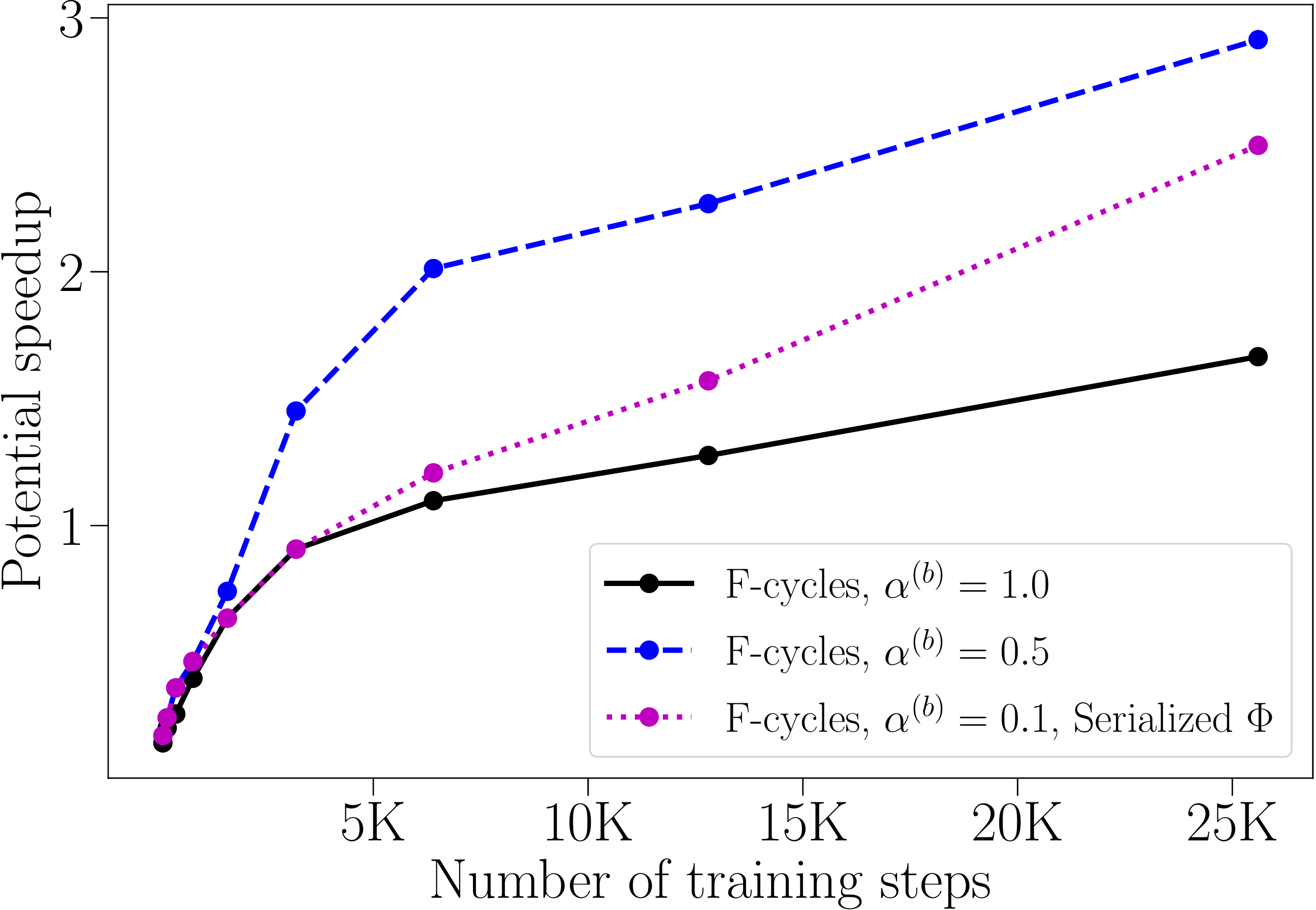} 
   \caption{Potential speedup for three-layer model problem, using Solver
   1 from Section \ref{sec:incalpha} with data from Tables
   \ref{tab:FcycleIncreasingAlpha} and
   \ref{tab:FcycleIncreasingAlpha_SmallBaseAlpha}, and Solver 3 from Section
   \ref{sec:3layer-serialize} with the data from Table
   \ref{tab:FcycleSerializedIncreasingAlpha_SmallBaseAlpha}. }
   \label{fig:threelayer_speedup} 
\end{figure}

Figure \ref{fig:threelayer_speedup} depicts the potential speedup for Solver 1 from Tables
\ref{tab:FcycleIncreasingAlpha}, \ref{tab:FcycleIncreasingAlpha_SmallBaseAlpha}
and Solver 3 from Table \ref{tab:FcycleSerializedIncreasingAlpha_SmallBaseAlpha}.
Although this is a simple model problem, these results indicate that the
observed convergence behavior offers promise for a useful speedup on more
real-world examples.  Also, the trend lines indicate that larger problem sizes
would offer improved potential speedups.

\begin{remark} \label{rem:future_speedup} We note that these speedup estimates
   are promising, there is much room for improvement.  This paper represents an
   initial work meant to show feasibility of the approach.  Future research
   will focus on improving the potential speedup, by for instance, considering
   coarsening the network on coarse levels, where the accuracy requirements of
   the training steps are reduced.  Such approaches have lead to significant
   speedup improvements for PDEs, and are expected to do so here.  Moreover, 
   there is existing work \cite{ChGoSh2015} on transferring a network's content
   between different sized networks.  
\end{remark}

\begin{remark} The computed potential speedup here assumes that the problem has
   access to as much parallelism as needed.  For $\N{0}$ training steps on
   level $0$, that means having $\N{0}/\cf{0}$ processors available.  If fewer
   processors are available, then the potential speedup would be reduced.
\end{remark}

\subsection{Larger Example}

In this section, we will investigate a larger example using the solvers
developed above for the small, three-layer model problem.

In particular, we will examine a four-layer network with $\globalwsize =
12\,032$ total weights.  The training data is a synthetic dataset for binary
addition created using a random number generator.  The network takes two random
binary numbers (12 bits in length) and then computes the 12 bit binary sum.
The error produced by the network is the difference between the computed and
actual sum.  The network is fully connected with an input layer of $\netsize{0}
= 24$ nodes (12 for each binary number), first hidden layer of $\netsize{1} =
128$ nodes, second hidden layer of $\netsize{2} = 64$ nodes and an output layer
of $\netsize{3} = 12$ nodes.  Note, $24\cdot128 + 128\cdot64 + 64\cdot12 =
12,032$.  This example is based on Section 4 from \cite{Trask_4Layer}, where the code of the
serial training algorithm is as the ``Plain Vanilla" algorithm based on
backpropagation.  The same sigmoid threshold function as in
Section \ref{sec:ann_intro} is used.  

We begin our experiments with the Naive Solver, using the $\basealpha=0.1$ from
\cite{Trask_4Layer}.  Remember for the Naive Solver, $\zl{\ell}{i}$ is defined
to contain all the training data sets, which here means $\xl{\ell}{i}$ is
defined as 500 pairs of random 12 bit numbers and $\yl{\ell}{i}$ then
corresponds to the sum of those 500 pairs.  

Table \ref{tab:2levelnaive_LargeBaseAlpha} depicts the two-level results for
the Naive Solver.  The results show some promise, but also large unscalable
growth in iteration counts.
\begin{table}[!ht]
  \centering
  \begin{tabular}{c | cccccc }
  \toprule 
    $\N{0}$       &  40  &  80  & 160  & 320  & 640  & 1280 \\
  \midrule       
    Iters       &  9   &  13  &  19  &  36  &  48  &  50  \\
    $\rho$      & 0.14 & 0.27 & 0.38 & 0.63 & 0.69 & 0.88 \\
  \bottomrule
  \end{tabular}
   \caption{Four-layer model problem, Naive solver, Two-level results, $\basealpha = 0.1$}
   \label{tab:2levelnaive_LargeBaseAlpha}
\end{table}

Following the direction of Section \ref{sec:incalpha}, we apply Solver 1,
which increases $\alphal{\ell}$ on coarse-levels.
We continue to use $\basealpha=0.1$ and set $\maxalpha=0.2$ (based on
experiments).  Table \ref{tab:2level_IncreasingAlpha_LargeBaseAlpha} depicts these two-level
Solver 1 results, which indicate some stability issues
for $\Phil{1}{i}$, because the convergence rate is over 1.0.
\begin{table}[!ht]
  \centering
  \begin{tabular}{c | cccccc }
  \toprule 
    $\N{0}$       &  40  &  80  & 160  & 320  & 640  & 1280 \\
  \midrule       
    Iters       &  9   &  50  &  50  &  50  &  50  &  50  \\
    $\rho$      & 0.13 & 1.01 & 1.02 & 1.02 & 1.02 & 1.02 \\
  \bottomrule
  \end{tabular}
   \caption{Four-layer model problem, Solver 1, Two-level results, $\basealpha = 0.1$}
   \label{tab:2level_IncreasingAlpha_LargeBaseAlpha}
\end{table}

Thus, we experiment with a smaller $\basealpha$.  Taking the Naive Solver, but
changing $\basealpha$ to $0.025$, yields the results depicted in Table
\ref{tab:2levelnaive_SmallBaseAlpha}.
\begin{table}[!ht]
  \centering
  \begin{tabular}{c | cccccc }
  \toprule 
    $\N{0}$       &  40  &  80  & 160  & 320  & 640  & 1280 \\
  \midrule       
    Iters       &  8   &  12  &  17  &  23  &  31  &  43  \\
    $\rho$      & 0.11 & 0.24 & 0.37 & 0.49 & 0.61 & 0.70 \\
  \bottomrule
  \end{tabular}
  \caption{Four-layer model problem, Naive Solver, Two-level results, $\basealpha = 0.025$}
  \label{tab:2levelnaive_SmallBaseAlpha}
\end{table}

Again following the direction of Section \ref{sec:incalpha}, we take this
smaller $\basealpha=0.025$ and experiment with it in the context of Solver 1,
continuing to use $\maxalpha=0.2$.  Table
\ref{tab:2level_IncreasingAlpha_SmallBaseAlpha} depicts these results, but
again, we see unscalable iteration growth.  There is a difficulty here for
MGRIT that was not present for the three-layer example.
\begin{table}[!ht]
  \centering
  \begin{tabular}{c | cccccc }
  \toprule 
    $\N{0}$       &  40  &  80  & 160  & 320  & 640  & 1280 \\
  \midrule       
    Iters       &  5   &  7   &  11  &  16  &  25  &  37  \\
    $\rho$      & 0.06 & 0.12 & 0.23 & 0.36 & 0.53 & 0.63 \\
  \bottomrule
  \end{tabular}
  \caption{Four-layer model problem, Solver 1, Two-level results, $\basealpha = 0.025$}
  \label{tab:2level_IncreasingAlpha_SmallBaseAlpha}
\end{table}

\subsection{A Better Solver}
\label{sec:better_solver}

To design a better solver, we begin by simply trying to understand which
problem trait(s) make the performance of MGRIT differ to that observed for
the small three-layer problem.  

One key difference is the number of training examples.  Here, the network is
trained on 500 examples, while before for the three-layer example, only 4
examples were used.  Thus, we take Solver 1 with $\basealpha=0.025$ and
$\maxalpha=0.2$, but cut $\xl{\ell}{i}$ to contain only 50 training examples.
Using this modified Solver 1 we get the results depicted in Table
\ref{tab:2level_IncreasingAlpha_SmallBaseAlpha_50examples}.  The reduction in
training examples has restored scalable iteration counts.
\begin{table}[!ht]
  \centering
  \begin{tabular}{c | cccccccc }
  \toprule 
    $\N{0}$       & 100  & 200  & 400  & 800  & 1600 & 3200 & 6400 & 12800 \\
  \midrule       
    Iters       &  4   &  5   &  7   &  11  &  14  &  13  &  14  &   14  \\
    $\rho$      & 0.04 & 0.07 & 0.15 & 0.27 & 0.35 & 0.37 & 0.39 &  0.39 \\
  \bottomrule
  \end{tabular}
  \caption{Four-layer model problem, Solver 1, Two-level results, 50 training examples,
   $\basealpha = 0.025$}
  \label{tab:2level_IncreasingAlpha_SmallBaseAlpha_50examples}
\end{table}

The use of 50 training set examples is not an attractive approach, because it
limits the ability to train a network.  However, it does indicate potential
promise for serializing over the training set data as in Section
\ref{sec:3layer-serialize}.  We explore this idea in the next section. 

\subsubsection{Serialized Training Sets}
\label{sec:serialized-4layer}

In this section, we experiment with Solver 2, which serializes over the
training data, and continue to use $\basealpha = 0.025$ and $\maxalpha = 0.2$.
In other words, the first 500 training steps will train over the first 500
training data sets; training step 501 will loop back to the first training data
set; and so on.  Because of the smaller training data set size per step, we are
able to run out to larger problem sizes.

Table \ref{tab:2level_Serialized_IncreasingAlpha_SmallBaseAlpha} depicts two-level 
results for Solver 2.  The iteration counts are good, but not quite
scalable, showing some growth. 
\begin{table}[!ht]
  \centering
  \begin{tabular}{c | cccccccc }
  \toprule 
     $\N{0}$       & 100  & 200  & 400  & 800  & 1600 & 3200 & 6400  & 12800  \\
  \midrule       
     Iters       &  4   &  4   &  5   &  5   &  6   &  6   &  7    & 11    \\
     $\rho$      & 0.05 & 0.07 & 0.07 & 0.09 & 0.10 & 0.13 & 0.18  & 0.29   \\
  \bottomrule
  \end{tabular}
  \caption{Four-layer model problem, Solver 2, Two-level results, $\basealpha=0.025$}
  \label{tab:2level_Serialized_IncreasingAlpha_SmallBaseAlpha}
\end{table}

Table \ref{tab:Fcycle_Serialized_IncreasingAlpha_SmallBaseAlpha} depicts
F-cycle results for Solver 2.  They are also not scalable, showing a slow
growth, but the iteration counts are modest enough to potentially be useful.
Table \ref{tab:Vcycle_Serialized_IncreasingAlpha_SmallBaseAlpha} depicts
results for Solver 2 for V-cycles.  The deterioration for larger problem sizes
and F-cycles in Table
\ref{tab:Fcycle_Serialized_IncreasingAlpha_SmallBaseAlpha} is visible when
compared to the two-level results in Table
\ref{tab:2level_Serialized_IncreasingAlpha_SmallBaseAlpha}.  To investigate
possible remedies, we again explore the possibility of allowing $\alphal{\ell}$
to grow on more levels.  Table
\ref{tab:Fcycle_Serialized_IncreasingAlpha_SmallBaseAlpha_SlowAlphaInc} carries
out such a experiment, except this time we compute $\alphal{\ell}$ with
\begin{equation} \label{eqn:inc_alpha2} \alphal{\ell} = \min(1.25^{\ell}
\basealpha, \maxalpha), \end{equation} as opposed to equation
(\ref{eqn:inc_alpha}), which uses powers of 2.  This allows $\alphal{\ell}$ to
grow on all levels for all experiments. The F-cycle results now show no
deterioration from the two-level results in Table
\ref{tab:2level_Serialized_IncreasingAlpha_SmallBaseAlpha} for the large
problem sizes.  Unfortunately, this choice of $\alphal{\ell}$ did not show a
significant improvement for V-cycles.  We leave a further investigation of
optimal coarse level $\alphal{\ell}$ as future work, because a deeper
theoretical understanding is required to meaningfully guide our tests.  
\begin{table}[!ht]
  \centering
  \begin{tabular}{c | cccccccc }
  \toprule 
    $\N{0}$     & 100  & 200  & 400  & 800  & 1600 & 3200 & 6400 & 12800 \\
  \midrule                                                               
    Iters       &  4   &  4   &  5   &  5   &  6   &  6   &  8   & 13    \\
    $\rho$      & 0.05 & 0.07 & 0.07 & 0.09 & 0.10 & 0.13 & 0.17 & 0.34  \\
  \bottomrule
  \end{tabular}
  \caption{Four-layer model problem, Solver 2, F-cycle results, $\basealpha=0.025$}
  \label{tab:Fcycle_Serialized_IncreasingAlpha_SmallBaseAlpha}
\end{table}
\begin{table}[!ht]
  \centering
  \begin{tabular}{c | cccccccc }
  \toprule 
    $\N{0}$     & 100  & 200  & 400  & 800  & 1600 & 3200 & 6400 & 12800 \\
  \midrule                                                               
    Iters       &  4   &  5   &  6   &  7   &  8   &  9   &  11  & 20    \\
    $\rho$      & 0.05 & 0.08 & 0.10 & 0.13 & 0.17 & 0.20 & 0.28 & 0.47  \\
  \bottomrule
  \end{tabular}
  \caption{Four-layer model problem, Solver 2, V-cycle results, $\basealpha=0.025$}
  \label{tab:Vcycle_Serialized_IncreasingAlpha_SmallBaseAlpha}
\end{table}
\begin{table}[!ht]
  \centering
  \begin{tabular}{c | cccccccc }
  \toprule 
    $\N{0}$     & 100  & 200  & 400  & 800  & 1600 & 3200 & 6400 & 12800 \\
  \midrule      
    Iters       &  4   &  4   &  4   &  5   &  6   &  6   &  7   &   10  \\
    $\rho$      & 0.04 & 0.06 & 0.07 & 0.08 & 0.09 & 0.12 & 0.14 &  0.27 \\
  \bottomrule
  \end{tabular}
   \caption{Four-layer model problem, Solver 2, F-cycle results, $\basealpha=0.025$, Slow $\alphal{\ell}$ increase}
  \label{tab:Fcycle_Serialized_IncreasingAlpha_SmallBaseAlpha_SlowAlphaInc}
\end{table}

\begin{remark} Further reducing the factor of increase from 1.25 in equation
   (\ref{eqn:inc_alpha2}) is not an option.  If we revert to the naive strategy
   of fixing $\alphal{\ell}$ to be $\basealpha$ on all levels for Solver 2,
   even the two-level solver shows unscalable iteration counts.  While we omit
these results for brevity, we note that the key insight of increasing
$\alphal{\ell}$ on coarse levels continues to be critical.  \end{remark}

\subsection{Potential Speedup}

\begin{figure}
  \center
  \includegraphics[width = 0.6\textwidth]{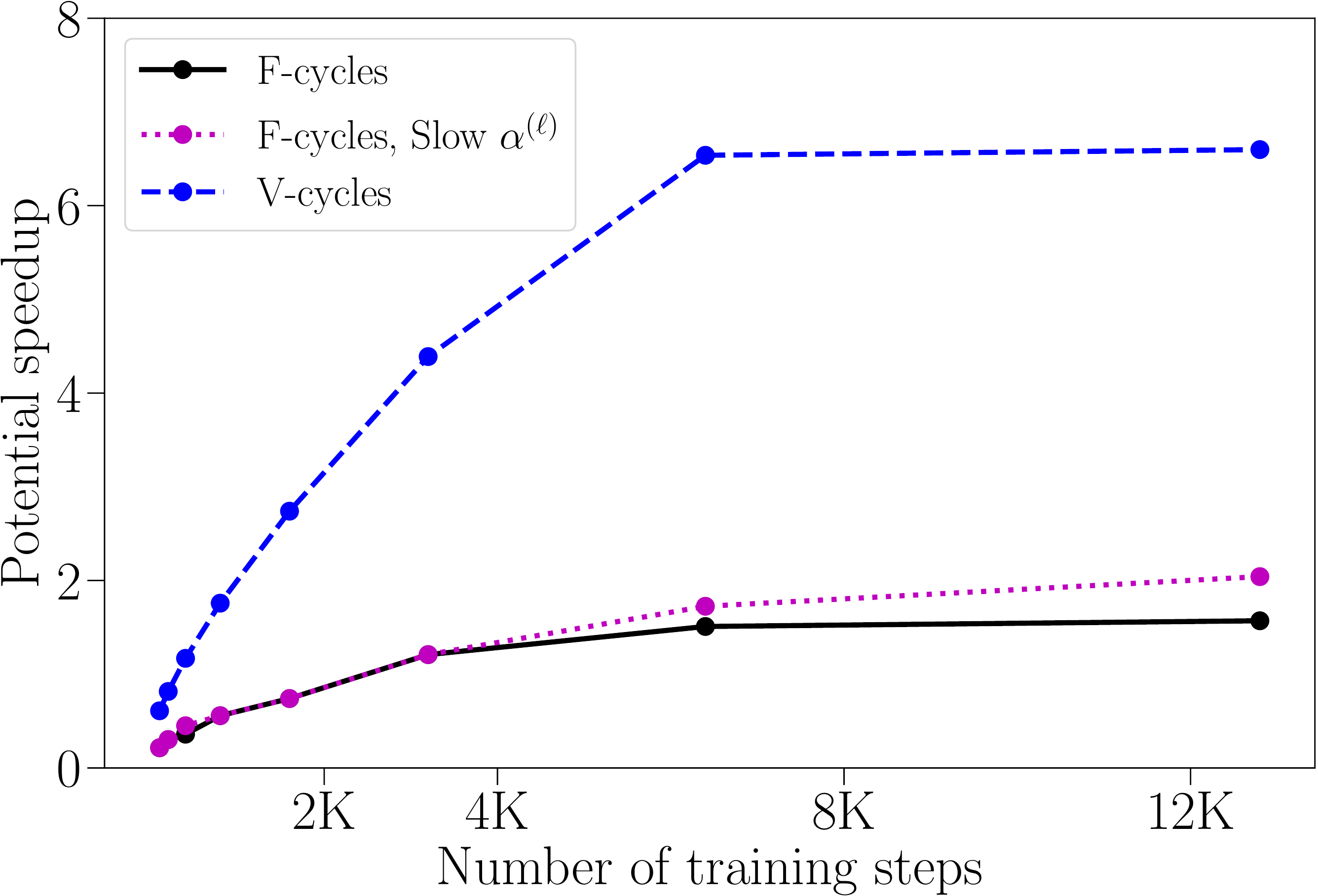}
  \caption{Potential speedup for four-layer model problem, using 
   Solver 2 and the data 
   from Tables \ref{tab:Fcycle_Serialized_IncreasingAlpha_SmallBaseAlpha},
   \ref{tab:Vcycle_Serialized_IncreasingAlpha_SmallBaseAlpha},
   and \ref{tab:Fcycle_Serialized_IncreasingAlpha_SmallBaseAlpha_SlowAlphaInc}. } 
  \label{fig:fourlayer_speedup}
\end{figure}

Figure \ref{fig:fourlayer_speedup} depicts the potential speedup for Solver 2 using 
Tables \ref{tab:Vcycle_Serialized_IncreasingAlpha_SmallBaseAlpha} and
\ref{tab:Fcycle_Serialized_IncreasingAlpha_SmallBaseAlpha}, respectively.  The
addition of another layer and many more weights to the network has not 
significantly impacted the possibility for speedup when compared 
to the three-layer example (see Section \ref{sec:3layer_speedup}).
We note that there is an unexplained deterioration in V-cycle convergence
for the final data point. Also, the trend lines indicate that larger problem sizes
would offer improved potential speedups.

\section{Conclusion}
\label{sec:conclusion}

We have shown a method that is able to parallelize over training steps (i.e.,
training runs) for artificial neural networks and also converge to the exact
same solution as that found by traditional sequential training algorithms.  To
the best of our knowledge, this is the first successful application of such
parallel-in-time techniques to neural networks. 

The method relies on a multigrid reduction technique to construct a hierarchy
of levels, where each level contains a fraction of the training steps found on
the preceding level.  The coarse levels serve to accelerate the solution on the
finest level.  The strategies developed here will hopefully lay the groundwork
for application to more real-world applications later.

The chief challenge faced is the definition of the coarse-level propagator
$\Phi$.  Proposed solutions for constructing such coarse-level $\Phi$ are
treating the learning rate in a way similar to time-step sizes and using double
applications of $\Phi$ to provide for more effective MGRIT relaxation.  We also
explored serializing over all training examples, to allow for more potential
parallelism.  Solver 3 was the most robust, and the results indicate that MGRIT
is able to offer a large possible speedup for the considered model problems.

Future research will focus on examining more challenging neural networks and on
improving the potential speedup (See Remark \ref{rem:future_speedup}).  Another
path for improving speedup is the development of a deeper theoretical
understanding of this approach.  This will help, for example, in choosing
$\alphal{\ell}$ and designing $\Phill{\ell}$.  Here, we will look to leverage
the existing MGRIT theory in \cite{Do2016}.


\bibliographystyle{siam}
\bibliography{refs}

\newpage
\appendix

\section{Code for Three-Layer Problem}
\label{sec:3layer_code}
As an aide to the reader, we reproduce the serial training algorithm (taken
directly from \cite{Trask_3Layer}) that we base the three-layer model problem
on from Section \ref{sec:ann_intro}.  The code is in Python.

\begin{lstlisting}
import numpy as np

def nonlin(x,deriv=False):
   if(deriv==True):
       return x*(1-x)

   return 1/(1+np.exp(-x))
    
X = np.array([[0,0,1],
              [0,1,1],
              [1,0,1],
              [1,1,1]])
y = np.array([[0],
              [1],
              [1],
              [0]])

# randomly initialize our weights with mean 0
np.random.seed(1)
syn0 = 2*np.random.random((3,4)) - 1
syn1 = 2*np.random.random((4,1)) - 1

for j in xrange(60000):

   # Feed forward through layers 0, 1, and 2
    l0 = X
    l1 = nonlin(np.dot(l0,syn0))
    l2 = nonlin(np.dot(l1,syn1))

    # how much did we miss the target value?
    l2_error = y - l2
    
    if (j% 10000) == 0:
        print "Error:" + str(np.mean(np.abs(l2_error)))
        
    # in what direction is the target value? 
    l2_delta = l2_error*nonlin(l2,deriv=True)

    # how much did each l1 value contribute to the l2 error?
    l1_error = l2_delta.dot(syn1.T)
    
    # in what direction is the target l1?  
    l1_delta = l1_error * nonlin(l1,deriv=True)
    
    alpha = 1.0
    syn1 += alpha*l1.T.dot(l2_delta)
    syn0 += alpha*l0.T.dot(l1_delta)
\end{lstlisting}


\end{document}